\documentclass[12pt]{article}

\usepackage{amssymb}

\usepackage{amsmath}

\usepackage{amsthm}

\newcommand{\cA}{{\cal A}}

\newcommand{\cO}{{\cal O}}

\newcommand{\cL}{{\cal L}}

\newcommand{\cN}{{\cal N}}

\newcommand{\cF}{{\cal F}}

\newcommand{\cQ}{{\cal Q}}
\newcommand{\cR}{{\cal R}}
\newcommand{\cS}{{\cal S}}

\newcommand{\cT}{{\cal T}}

\newcommand{\cV}{{\cal V}}

\newcommand{\cW}{{\cal W}}

\newcommand{\cX}{{\cal X}}

\newcommand{\cY}{{\cal Y}}

\newcommand{\cZ}{{\cal Z}}

\renewcommand{\AA}{{\mathbb A}}

\newcommand{\ZZ}{{\mathbb Z}}

\newcommand{\CC}{{\mathbb C}}

\newcommand{\on}{\operatorname}

\newcommand{\Qlb}{\mathbb{\bar Q}_\ell}

\newcommand{\Gm}{\mathbb{G}_m}

\newcommand{\toup}[1]{\stackrel{#1}{\to}}

\newcommand{\hook}[1]{\stackrel{#1}{\hookrightarrow}}

\newcommand{\Sp}{\on{\mathbb{S}p}}

\newcommand{\Hom}{\on{Hom}}

\newcommand{\Ext}{\on{Ext}}

\newcommand{\Sym}{\on{Sym}}

\newcommand{\SO}{\on{S\mathbb{O}}}

\newcommand{\Ker}{\on{Ker}}

\newcommand{\Aut}{\on{Aut}}

\newcommand{\RG}{\on{R\Gamma}}

\newcommand{\Bun}{\on{Bun}}

\newcommand{\Bunt}{\on{\widetilde\Bun}}

\newcommand{\Spec}{\on{Spec}}

\newcommand{\QED}{$\square$}

\newcommand{\Fp}{\mathbb{F}_p}  

\newcommand{\iso}{{\widetilde\to}}

\newcommand{\comp}{\circ}

\renewcommand{\H}{{\on{H}}}   

\newcommand{\R}{\on{R}\!}   

\newcommand{\D}{\on{D}}       

\newcommand{\select}[1]{{\it{#1}}}

\renewcommand{\div}{\on{div}}

\newcommand{\<}{\langle}

\renewcommand{\>}{\rangle}

\newcommand{\ev}{\mathit{ev}}

\newcommand{\dimrel}{\on{dim.rel}}

\newcommand{\Loc}{\mathrm{Loc}}

\newcommand{\GL}{\on{\mathbb{G}L}}

\newtheorem{Lm}{Lemma}

\newtheorem{Th}{Theorem}

\newtheorem{Pp}{Proposition}

\theoremstyle{remark}

\newtheorem{Rem}{Remark}

\theoremstyle{definition}

\newtheorem{Def}{Definition}

\newenvironment{Prf}{\par\noindent {\it Proof }}{\QED}

\def\det{\mathrm{det}}
\newcommand{\s}[1]{\langle #1 \rangle}
\newcommand{\vag}[1]{\overset{#1}{\leftarrow}}
\newcommand{\vad}[1]{\overset{#1}{\rightarrow}}
\def\Sym{\mathrm{Sym}}

\def\Whit{\mathrm{Whit}}

\begin{document}

\title{Compatibility of the Theta correspondence with the Whittaker functors}
\author{Vincent Lafforgue and Sergey Lysenko}
\date{}
\maketitle

We prove in this note that the global geometric theta lifting 
for the pair $(H, G)$ is compatible with the Whittaker normalization, where $(H,G)=(\SO_{2n},\Sp_{2n})$,  $(\Sp_{2n},\SO_{2n+2})$, or $(\GL_{n},\GL_{n+1})$. More precisely, let $k$ be an algebraically closed field of characteristic $p>2$.
Let $X$ be a smooth projective connected curve over $k$.
For a stack $S$ write $\D(S)$ for the derived category of \'etale constructible $\Qlb$-sheaves on $S$. For a reductive group $G$ over $k$ write $\Bun_G$ for the stack of $G$-torsors on $X$. The usual Whittaker distribution admits a natural geometrization $\Whit_{G}:\D(\Bun_{G})\to \D(\Spec k)$. 

 We construct an isomorphism of functors between $\Whit_{G}\circ F$ and $\Whit_{H}$ where $F:\D(\Bun_{H}) \to \D(\Bun_{G})$ is the theta lifting functor (cf. Theorems 1, 2 and 3).

This result at the level of functions (on $\Bun_{H}(k)$ and $\Bun_{G}(k)$ when $k$ is a finite field) is well known since a long time and the geometrization of the argument is straightforward.  We  wrote this note   for the following reason.

Our proof holds also for $k=\CC$ in the setting of $D$-modules. In this case for a reductive group $G$, Beilinson and Drinfeld proposed a conjecture, which (in a form that should be made more precise) says that there exists an equivalence $\alpha_{G}$ between  the derived category of $\D$-modules on $\Bun_{G}$ and  the derived category of $\cO$-modules on $\Loc_{\check{G}}$. Here $\Loc_{\check{G}}$ is the stack of $\check{G}$-local systems on $X$, and $\check{G}$ is the Langlands dual group to $G$. Moreover, $\Whit_{G}$ should be the composition $\D(D\mathrm{-mod}(\Bun_{G}))\vad{\alpha_{G}} \D(\Loc_{\check{G}}, \cO)\vad{\R\Gamma} \D(\Spec \CC)$. 

A morphism $\gamma: \check{H}\to \check{G}$ gives rise to the extension of scalars morphism $\bar\gamma: \Loc_{\check{H}}\to \Loc_{\check{G}}$. The functor $\bar\gamma_*:\D(\Loc_{\check{H}},\cO)\to \D(\Loc_{\check{G}}, \cO)$ should give rise to the Langlands functoriality functor 
$$
\gamma_L=\alpha_G^{-1}\comp \bar\gamma_*\comp \alpha_H: \D(D\mathrm{-mod}(\Bun_H))\to\D(D\mathrm{-mod}(\Bun_G))
$$ 
compatible with the action of Hecke functors. 

  In the cases $(H,G)=(\SO_{2n},\Sp_{2n})$,  $(\Sp_{2n},\SO_{2n+2})$ or $(\GL_{n},\GL_{n+1})$ the compatibility of the theta lifting functor $F : \D(D\mathrm{-mod}(\Bun_{H})) \to D(D\mathrm{-mod}(\Bun_{G}))$ with the Hecke functors (\cite{L1}) and the compatibility of $F$ with the Whittaker functors (proved in this paper) indicate that $F$ should be the Langlands functoriality functor. 
  
{\scshape Notation.} From now on $k$ denotes an algebraically closed field of characteristic $p>2$, all the stacks we consider are defined over $k$. Let $X$ be a smooth projective curve of genus $g$. Fix a prime $\ell\neq p$ and a non-trivial character $\psi : \Fp\to \Qlb^{*}$, and denote by $\cL_{\psi}$ the corresponding Artin-Schreier sheaf on $\AA^{1}$. Since $k$ is algebraically closed, we systematically ignore the Tate twists. 

 For a $k$-stack locally of finite type $S$ write simply $\D(S)$ for the category introduced in (\cite{LO}, Remark~3.21) and denoted $\D_c(S,\Qlb)$ in \select{loc.cit}. It should be thought of as the unbounded derived category of constructible $\Qlb$-sheaves on $S$. For $\ast=+,-, b$ we have the full triangulated subcategory $\D^{\ast}(S)\subset \D(S)$ denoted $\D_c^{\ast}(S,\Qlb)$ in \select{loc.cit.} Write 
$\D^{\ast}(S)_!\subset \D^{\ast}(S)$ for the full subcategory of objects which are extensions by zero from some open substack of finite type. Write $\D^{\prec}(S)\subset \D(S)$ for the full subcategory of complexes $K\in D(S)$ such that for any open substack $U\subset S$ of finite type we have $K\mid_U\in D^-(U)$. 

For any vector space (or bundle) $E$, we define $\Sym^{2}(E)$ and $\Lambda^{2}(E)$ as quotients of $E\otimes E$ (and denote by 
$x.y$ and $x\wedge y$ the images of $x\otimes y$)  
and we will use in this article the embeddings 
\begin{eqnarray}
\label{plongement}
\begin{matrix}
\Sym^{2}(E) &\to & E\otimes E &\ \  \text{and}\ \  &  \Lambda^{2}(E)&\to & E\otimes E \\
x.y&\mapsto & \frac{x\otimes y+y\otimes x}{2} && x \wedge y&\mapsto & \frac{x\otimes y-y\otimes x}{2}
\end{matrix}
\end{eqnarray}
        
\section{Whittaker functors}

Let $G$ be a reductive group over $k$. We pick a maximal torus and a Borel subgroup $T\subset B\subset G$ and  we denote by $\Delta_{G}$ the set of simple roots of $G$. 
The Whittaker functor 
$$
\Whit_{G}:\D^{\prec}(\Bun_{G})\to \D^-(\Spec k)
$$ 
is defined as follows. Write $\Omega$ for the canonical line bundle on $X$. Pick a $T$-torsor $\cF_T$ on $X$ with a trivial conductor, that is, for each  $\check \alpha\in \Delta_{G}$ it is equipped with an isomorphism $\delta_{\check \alpha}:\cL^{\check{\alpha}}_{\cF_T}\to \Omega$. Here $\cL^{\check{\alpha}}_{\cF_T}$ is the line bundle obtained from $\cF_T$ via extension of scalars $T\toup{\check{\alpha}}\Gm$. Let $\Bun_{N}^{\cF_{T}}$ be the stack classifying a $B$-torsor $\cF_B$ together with an isomorphism 
$$
\zeta : \cF_B\times_B T\,\iso\,\cF_T
$$ 
Let $\epsilon:\Bun_{N}^{\cF_{T}}\to \mathbb{A}^{1}$ be the evaluation map (cf. \cite{whittaker}, 4.3.1 where it is denoted $\ev_{\tilde\omega}$). 
Just recall that 
for each $\check{\alpha}\in \Delta_{G}$ the class of the extension 
of $\cO$ by $\Omega$ associated to $\cF_{B}$, $\zeta$ and $\delta_{\check{\alpha}}$ 
 gives $\epsilon_{\check{\alpha}}:\Bun_{N}^{\cF_{T}}\to \mathbb{A}^{1}$ and that $\epsilon=\sum_{\check{\alpha}\in \Delta_{G}}\epsilon_{\check{\alpha}}$. Write $\pi:\Bun_{N}^{\cF_{T}}\to \Bun_{G}$ for the extension of scalars 
 $(\cF_{B},\zeta)\mapsto \cF_B\times_B G$. 
 Set $P_{\psi}^{0}=\epsilon^*\cL_{\psi}[d_N]$, where $d_N=\dim\Bun_N^{\cF_T}$. Let $d_G=\dim\Bun_G$. As in (\cite{sergey}, Definition~2) for $\cF\in \D^{\prec}(\Bun_G)$ set
\begin{eqnarray}\label{def-whit}
\Whit_{G}(\cF)=\mathrm{R}\Gamma_{c}(\Bun_{N}^{\cF_{T}},P_{\psi}^{0}\otimes \pi^{*}(\cF))[-d_G]
\end{eqnarray}

\begin{Rem} The collection $(\cF_{T},(\delta_{\check \alpha})_{\check \alpha\in \Delta_{G}})$ as above exists, because $k$ is algebraically closed, and one can take $\cF_T=(\sqrt{\Omega})^{2\rho}$ for some square root $\sqrt{\Omega}$ of $\Omega$.
One has an exact sequence of abelian group schemes $1\to Z  \to T\vad{\prod \check{\alpha}} \Gm^{\Delta_{G}}\to 1$ where $Z$ denotes the center of $G$. 
So, two choices of the collection $(\cF_{T},(\delta_{\check \alpha})_{\check \alpha\in \Delta_{G}})$ are related by a point of $\Bun_{Z}(k)$ and the associated Whittaker functors  are isomorphic up to the automophism of $\Bun_{G}$ given by tensoring with the corresponding $Z$-torsor. 
\end{Rem} 
\begin{Rem} \label{remarque2} 
When $\cF_{T}$ is fixed, the functor $\Whit_{G}:\D^{\prec}(\Bun_{G})\to \D^-(\Spec k)$ does not depend, up to isomorphism, on the choice of the isomorphisms  $(\delta_{\check \alpha})_{\check \alpha\in \Delta_{G}}$. 
That is, for any $(\lambda_{\check \alpha})_{\check \alpha\in \Delta_{G}}\in (k^{*})^{\Delta_{G}}$, the functors associated to $(\cF_{T},(\delta_{\check \alpha})_{\check \alpha\in \Delta_{G}})$ and  $(\cF_{T},(\lambda_{\check \alpha}\delta_{\check \alpha})_{\check \alpha\in \Delta_{G}})$ are isomorphic. Indeed, the two diagrams 
$\Bun_{G}\vag{\pi}\Bun_{N}^{\cF_{T}}\vad{\epsilon} \mathbb{A}^{1}$
associated to $(\delta_{\check \alpha})_{\check \alpha\in \Delta_{G}}$ and $
(\lambda_{\check \alpha}\delta_{\check \alpha})_{\check \alpha\in \Delta_{G}}$ are isomorphic for the following reason. Since 
$k$ is algebraically closed, $T(k)\to (k^{*})^{\Delta_{G}}$ is surjective. We pick any preimage $\gamma\in T(k)$ of $(\lambda_{\check \alpha})_{\check \alpha\in \Delta_{G}}$ and get the automorphism 
$(\cF_{B},\zeta)\mapsto (\cF_{B},\gamma\zeta)$ 
of $\Bun_{N}^{\cF_{T}}$, which together with the idendity of $\Bun_{G}$ and $\mathbb{A}^{1}$ intertwines the two diagrams. 
\end{Rem}

\subsection{Whittaker functor for $\GL_{n}$}

 For $i,j\in \ZZ$ with $i\leq j$ we denote by $\cN_{i,j}$ the stack classifying the extensions of 
 $\Omega^{i}$ by $\Omega^{i+1}$  ... by  $\Omega^{j}$, i.e. classifying 
 a vector bundle $E_{j-i+1}$ on $X$ with a complete flag of vector subbundles $0=E_{0}\subset E_{1}\subset ...\subset E_{j-i+1}$ together with isomorphisms
$E_{k+1}/E_{k}\simeq \Omega^{j-k}$ for $k=0,...,j-i$. Write $\epsilon_{i,j} :  \cN_{i,j}\to \AA^1$ for the map given by the sum of the classes in $\Ext^1(\cO, \Omega)\,\iso\,\AA^1$ of the extensions $0\to E_{k+1}/E_k\to E_{k+2}/E_k\to E_{k+2}/E_{k+1}\to 0$ for $k=0,\ldots, j-i-1$.

For $G=\GL_{n}$, 
we consider the diagram $\Bun_{n}\vag{\pi_{0,n-1}} \cN_{0,n-1}\vad{\epsilon_{0,n-1}}\AA^{1}$, where $\pi_{0,n-1}:\cN_{0,n-1}\to \Bun_n$ is 
$(0=E_{0}\subset \dots \subset E_{n})\mapsto E_{n}$. This diagram is isomorphic to the diagram 
$\Bun_{G}\vag{\pi}\Bun_{N}^{\cF_{T}}\vad{\epsilon} \mathbb{A}^{1}$
associated to the  choice of $\cF_{T}$ whose image in $\Bun_{n}$ is $\Omega^{n-1}\oplus\Omega^{n-2}\oplus\ldots\oplus\cO$. 

Therefore the  functor 
$\Whit_{\GL_n}: \D^{\prec}(\Bun_{n})\to \D^-(\Spec k)$ associated to the above choice of $\cF_{T}$  is given by
$$
\Whit_{\GL_n}(\cF)=\mathrm{R}\Gamma_{c}(\cN_{0,n-1},\epsilon _{0,n-1}^{*}(\cL_{\psi})\otimes \pi_{0,n-1}^{*}(\cF))[\dim \cN_{0,n-1}-\dim\Bun_n]. 
$$

\begin{Rem} If $E$ is an irreducible rank $n$ local system on $X$ let $\Aut_E$ be the corresponding automorphic sheaf on $\Bun_n$ (cf. \cite{FGV-geom}) normalized to be perverse. Then $\Aut_E$ is equipped with a canonical isomorphism $\Whit_{\GL_{n}}(\Aut_E)\,\iso\, \Qlb$. This is
our motivation for the above shift normalization in~(\ref{def-whit}).
\end{Rem}

\subsection{Whittaker functor for $\Sp_{2n}$}

Write $G_n$ for the group scheme on $X$ of automorphisms of $\cO^n\oplus\Omega^n$ preserving the natural symplectic form $\wedge^2(\cO^n\oplus\Omega^n)\to\Omega$. The stack $\Bun_{G_n}$ of $G_n$-torsors on $X$ can be seen as the stack classifying vector bundles $M$ over $X$ of rank $2n$ equipped with a non-degenerate symplectic form $\Lambda^2 M\to \Omega$. 

The diagram $\Bun_{G_n}\vag{\pi_{G_n}}\cN_{G_n}\vad{\epsilon _{G_n}}\AA^{1}$ constructed in the next definition is isomorphic to the diagram 
$\Bun_{G}\vag{\pi}\Bun_{N}^{\cF_{T}}\vad{\epsilon} \mathbb{A}^{1}$
associated, for $G=G_{n}$, to the choice 
%
of $\cF_{T}$ whose image in $\Bun_{G_n}$ is $L\oplus L^*\otimes\Omega$ with $L=\Omega^n\oplus\Omega^{n-1}\oplus\ldots\oplus\Omega$
(with the natural  symplectic structure for which $L$ and $L^*\otimes\Omega$ are lagrangians).

\begin{Def} Let 
$\cN_{G_n}$ be the stack classifying $((L_{1},...,L_{n}),E)$, where 
$(0=L_{0}\subset L_{1}\subset ... \subset L_{n})\in  \cN_{1,n}$, and $E$ is an extension of $\cO_X$-modules
\begin{equation}
\label{ext_one}
0\to \Sym^{2} L_{n}\to E\to \Omega\to 0
\end{equation}

We associate to (\ref{ext_one}) an extension 
\begin{equation}
\label{ext_for_M}
0\to L_{n}\to M\to L_{n}^{*}\otimes \Omega
\to 0
\end{equation}
with  $M\in \Bun_{G_n}$ and $L_{n}$ lagrangian as follows. Equip $L_n\oplus L_n^*\otimes\Omega$ with the symplectic form $(l,l^*), (u,u^*)\mapsto \<l, u^*\>-\<u, l^*\>$ for $l,u\in L, l^*,u^*\in L^*$. Here $\<.,.\>$ is the canonical paring between $L_n$ and $L_n^*$. Using (\ref{plongement}), we consider (\ref{ext_one}) as a torsor on $X$ under the sheaf of symmetric morphisms $L_n^*\otimes\Omega\to L_n$. The latter sheaf acts naturally on $L_n\oplus L_n^*\otimes\Omega$ preserving the symplectic form. Then $M$ is the twisting of $L_n\oplus L_n^*\otimes\Omega$ by the above torsor. This defines a morphism  $\pi_{G_n} : \cN_{G_n}\to \Bun_{G_n}$. 

Note that the extension of $\Omega$ by $L_{n}\otimes L_{n}$ obtained from (\ref{ext_for_M}) is the push-forward of (\ref{ext_one}) by the embedding $\Sym^{2} L_{n}\to L_{n}\otimes L_{n}$ we have fixed in~(\ref{plongement}). 

Let $\epsilon _{G_n} : \cN_{G_n}\to \AA^{1}$ denote the sum of $\epsilon_{1,n}(L_{1},...,L_{n})$ with the class in $\Ext(\cO,\Omega)=\AA^{1}$ of the push-forward of (\ref{ext_one}) by $\Sym^{2 }L_{n}\to \Sym^{2}(L_{n}/L_{n-1})=\Omega^{2}$.
\end{Def}
The functor 
$\Whit_{G_n}: \D^{\prec}(\Bun_{G_n})\to \D^-(\Spec k)$ associated to the above choice of $\cF_{T}$  is given by
$$
\Whit_{G_n}(\cF)=\mathrm{R}\Gamma_{c}(\cN_{G_n},\epsilon _{G_n}^{*}(\cL_{\psi})\otimes \pi_{G_n}^{*}(\cF))[d_{N(G_n)}-d_{G_n}]
$$
with $d_{N(G_n)}=\dim \cN_{G_n}$ and $d_{G_n}=\dim\Bun_{G_n}$.

\subsection{Whittaker functor for $\SO_{2n}$ (first form)}

Let $H_n=\SO_{2n}$.The stack $\Bun_{H_n}$ of $H_n$-torsors can be seen as the stack classifying vector bundles $V$ over $X$ equipped with a non-degenerate symmetric form $\Sym^{2}V\to \cO$ and a compatible trivialization $\det V\,\iso\, \cO$. 

The diagram $\Bun_{H_n}\vag{\pi_{H_n}}\cN_{H_n}\vad{\epsilon _{H_n}}\AA^{1}$ constructed in the next definition is isomorphic to the diagram 
$\Bun_{G}\vag{\pi}\Bun_{N}^{\cF_{T}}\vad{\epsilon} \mathbb{A}^{1}$
associated, for $G=H_{n}$, to the choice 
of $\cF_{T}$ whose image in  $\Bun_{H_n}$ is $U\oplus U^*$ with $U=\Omega^{n-1}\oplus \Omega^{n-2}\oplus\dots \oplus \cO$ (with the natural symmetric structure for which $U$ and $ U^*$ are isotropic).

\begin{Def}\label{whit-so-1} Let 
$\cN_{H_n}$ be the stack classifying
$((U_{1},...,U_{n}),E)$, where $(U_{1},...,U_{n})\in \cN_{0,n-1}$ 
(i.e. we have a filtration 
$0=U_{0}\subset U_{1}\subset ... \subset U_{n}$ with $U_{i}/U_{i-1}\simeq \Omega^{n-i}$ for $i=1,...,n$), 
and $E$ is an extension of $\cO_X$-modules
\begin{equation}
\label{ext_for_SO}
0\to \Lambda^2 U_{n}\to E\to \cO\to 0
\end{equation} 

We associate to (\ref{ext_for_SO}) an extension 
\begin{equation}
\label{ext_for_V}
0\to U_{n}\to V\to U_{n}^{*}
\to 0
\end{equation}
with  $V\in \Bun_{H_n}$ and $U_{n}$ isotropic as follows. Equip $U_n\oplus U_n^*$ with the symmetric form given by $(u,u^*),(v,v^*)\mapsto \<u, v^*\>+\<v, u^*\>$ with $u,v\in U_n, u^*,v^*\in U_n^*$. Using (\ref{plongement}), we consider (\ref{ext_for_SO}) as a torsor under the sheaf of antisymmetric morphisms $U_n^*\to U_n$ of $\cO_X$-modules. This sheaf acts naturally on $U_n\oplus U_n^*$ preserving the symmetric form and the trivialization of $\det(U_n\oplus U_n^*)$. Then (\ref{ext_for_V}) is the twisting of 
$U_n\oplus U_n^*$ by the above torsor. This defines a morphism   $\pi_{H_n} : \cN_{H_n}\to \Bun_{H_n}$.

 Note that the extension of $\cO_X$ by $U_{n}\otimes U_{n}$ obtained from (\ref{ext_for_V})
is the push-forward of (\ref{ext_for_SO}) by the embedding $\Lambda^{2} U_{n}\to U_{n}\otimes U_{n}$ fixed in~(\ref{plongement}).

 For $\lambda\in k^*$ let $\epsilon _{H_n, \lambda} : \cN_{H_n}\to \AA^{1}$ be the sum of $\epsilon_{0,n-1}(U_{1},...,U_{n})$ with $\lambda u$, where $u\in \Ext(\cO,\Omega)=\AA^{1}$ is the class of the push-forward of (\ref{ext_for_SO}) by  $\Lambda^{2 }U_{n}\to \Lambda^{2}(U_{n}/U_{n-2})=\Omega$. Set $\epsilon _{H_n}=\epsilon_{H_n, 1}$.
\end{Def}
The functor $\Whit_{H_n}: \D^{\prec}(\Bun_{H_n})\to \D^-(\Spec k)$ associated to the above choice of $\cF_{T}$ sends $\cF\in D^{\prec}(\Bun_{H_n})$ to 
\begin{equation}
\label{Whit_H_n_def}
\Whit_{H_n}(\cF)=\mathrm{R}\Gamma_{c}(\cN_{H_n},\epsilon _{H_n}^{*}(\cL_{\psi})\otimes \pi_{H_n}^{*}(M))[d_{N(H_n)}-d_{H_n}]
\end{equation}
with $d_{N(H_n)}=\dim \cN_{H_n}$ and $d_{H_n}=\dim\Bun_{H_n}$. By Remark~\ref{remarque2}, if we replace in (\ref{Whit_H_n_def}) $\epsilon_{H_n}$ by $\epsilon_{H_n,\lambda}$ then the functor $\Whit_{H_n}$ gets replaced by an isomorphic one.

\subsection{Whittaker functor for $\SO_{2n}$ (second form)}

\begin{Def}\label{whit-so-2} Let 
$\widetilde\cN_{H_n}$ be the stack classifying
$(V_{1}\subset\ldots\subset V_{n}\subset V)$, where $V\in \Bun_{H_n}$, $V_{n}\subset V$ is a subbundle, $(V_{1},...,V_{n})\in  \cN_{0,n-1}$ (i.e. we have a filtration 
$0=V_{0}\subset V_{1}\subset ... \subset V_{n}$ with $V_{i}/V_{i-1}\simeq \Omega^{n-i}$ for $i=1,...,n$), 
and the composition 
$$
\Sym^{2}V_{n}\to \Sym^{2}V\to \cO
$$ 
coincides with $\Sym^{2}V_{n}\to \Sym^{2}(V_{n}/V_{n-1})=\cO$ (in particular $V_{n-1}$ is isotropic).  

The  morphism  $\widetilde\pi_{H_n} : \widetilde\cN_{H_n}\to \Bun_{H_n}$ sends $((V_{1},...,V_{n}),V)$ to $V$. The morphism
$\widetilde\epsilon _{H_n} : \widetilde\cN_{H_n}\to \AA^{1}$ is given by 
$\widetilde\epsilon _{H_n}((V_{1},...,V_{n}),V)=
\epsilon_{0,n-1}(V_{1},...,V_{n})$. 
\end{Def}

 Define a morphism $\kappa: \cN_{H_n}\to \widetilde\cN_{H_n}$ as follows.  
 Let $(U_{1},...,U_{n}),E)\in \cN_{H_n}$ and let $V$ be as in Definition~\ref{whit-so-1}. For $i=1,...,n-1$ define $V_i$
as the image of $U_{i}$ in $V$ and $V_{2n-i}$ as the orthogonal of $V_{i}$ in $V$. Then we have a filtration $$0=V_{0}\subset V_{1}\subset ...\subset V_{n-1}\subset V_{n+1}\subset ...\subset V_{2n-1}\subset V_{2n}=V.$$
Recall that we have an identification $U_{n}/U_{n-1}\simeq \cO$. 
The exact sequence $0\to U_n/U_{n-1}\to V_{n+1}/V_{n-1}\to V_{n+1}/U_n\to 0$ admits a unique splitting $s$ such that the image of $\cO=V_{n+1}/U_n\toup{s} V_{n+1}/V_{n-1}$ is isotropic. Thus, 
$V_{n+1}/V_{n-1}$ is canonically identified with $\cO\oplus \cO$ in such a way that the symmetric bilinear form 
$\Sym^{2}(\cO\oplus \cO)\to \cO$ becomes
$$
(1,0).(1,0)\mapsto 0,\ (1,0).(0,1)\mapsto 1,\ (0,1).(0,1)\mapsto 0 
$$
Under this identification $\cO=U_{n}/U_{n-1}\to V_{n+1}/V_{n-1}=\cO\oplus \cO$ sends $1$ to $(1,0)$. 

 Define $V_{n}$, equipped with $\cO\simeq V_{n}/V_{n-1}$ by the property that $\cO\simeq V_{n}/V_{n-1}\hook{} V_{n+1}/V_{n-1}$ sends $1$ to $(1,\frac{1}{2})\in \cO\oplus \cO$. The following is easy to check.
 
\begin{Lm}\label{lemme1-tilde} The map $\kappa 
: \cN_{H_n}\to \widetilde\cN_{H_n}$ is an isomorphism.
There exists $\lambda\in k^*$ such that $\widetilde\epsilon _{H_n}\circ \kappa
=\epsilon _{H_n, \lambda}$ and 
$\widetilde\pi _{H_n}\circ \kappa
=\pi _{H_n}$. \QED
\end{Lm}

 By Remark~\ref{remarque2}, if we replace in (\ref{Whit_H_n_def}) 
$\epsilon_{H_n}, \pi_{H_n}$ by $\tilde\epsilon_{H_n}, \tilde\pi_{H_n}$ then the functor $\Whit_{H_n}$ gets replaced by an isomorphic one.

\section{Main statements}

Write $\Bun_n$ for the stack of rank $n$ vector bundles on $X$. Let $\Bun_{P_n}$ be the stack classifying $L\in\Bun_n$ and an exact sequence $0\to \Sym^2 L\to ?\to\Omega\to 0$.
Remind the complex $S_{P,\psi}$ on $\Bun_{P_n}$ introduced in (\cite{sergey-theta}, 5.2). Let $\cV$ be the stack over $\Bun_n$ whose fibre over $L$ is $\Hom(L,\Omega)$. For $\cX_n=\cV\times_{\Bun_n}\Bun_{P_n}$ let 
 $p: \cX_n\to\Bun_{P_n}$ be the projection. Write $q: \cX_n\to\AA^1$ for the map sending $s\in\Hom(L,\Omega)$ to the pairing of $s\otimes s\in\Hom(\Sym^2L, \Omega^2)$ with the exact sequence $0\to \Sym^2 L\to ?\to\Omega\to 0$.
 Let $d_{\cX_n}$ be the ``corrected'' dimension of $\cX_n$, i.e. the locally constant function $\dim \Bun_{P_{n}}-\chi(L)$. 
  Set
$$
S_{P,\psi}=p_!q^*\cL_{\psi}[d_{\cX_n}].
$$

 Let $\cA$ be the line bundle on $\Bun_{G_n}$ whose fibre at $M$ is $\det\RG(X,M)$. Write $\Bunt_{G_n}$ for the gerb of square roots of $\cA$ and $\Aut$ for the theta-sheaf on $\Bunt_{G_n}$ (\cite{sergey-theta}, Definition~1). The projection $\nu_n: \Bun_{P_n}\to \Bun_{G_n}$ lifts naturally to a map $\tilde\nu_n: \Bun_{P_n}\to\Bunt_{G_n}$. In what follows, we pick an isomorphism\footnote{Once $\sqrt{-1}\in k$ is chosen, this isomorphism is well defined up to a sign.}
\begin{equation}
\label{iso_S_Ppsi}
S_{P,\psi}\,\iso\, \tilde\nu_n^*\Aut[\dimrel(\tilde\nu_n))]
\end{equation}
provided by (\cite{sergey-waldspurger}, Proposition~1). Here $\dimrel(\tilde\nu_n)$ is the relative dimension of $\tilde\nu_n$. The isomorphisms we construct below may depend on this choice.  
  
\subsection{From $\Sp_{2n}$ to $\SO_{2n+2}$}

Let $F:\D^-(\Bun_{G_n})_!\to \D^{\prec}(\Bun_{H_{n+1}})$ be the theta lifting functor introduced in (\cite{L1}, Definition~2).

\begin{Th}
\label{Th_1} The functors $\Whit_{H_{n+1}}\circ F$
and $\Whit_{G_n}$ from $\D^-(\Bun_{G_n})_!$ to $\D^-(\Spec k)$ are isomorphic. 
\end{Th}

Let $\cX$ be the stack classifying  $(M,(U_{1},...,U_{n+1}),E,s)$ with $M\in \Bun_{G_n}$, $(U_{1},...,U_{n+1})\in \cN_{0,n}$ (i.e. $U_{k+1}/U_{k}=\Omega^{n-k}$ for $k=0,...,n$), 
$E$ an extension $0\to \Lambda^{2}U_{n+1}\to E\to \cO\to 0$, and $s:U_{n+1}\to M$ a morphism of $\cO_X$-modules.

Let $\alpha_{\cX} : \cX\to \Bun_{G_n}$ be the morphism 
$(M,(U_{1},...,U_{n+1}),E,s)\mapsto M$. 
Let $\beta_{\cX} : \cX\to \AA^{1}$ be defined as follows. For 
$(M,(U_{1},...,U_{n+1}),E,s)\in \cX$, 
$$\beta_{\cX}(M,(U_{1},...,U_{n+1}),E,s)=\epsilon_{0,n}(U_{1},...,U_{n+1})+\gamma(E)-\s{E,\Lambda^{2}s}$$
where $\gamma(E)$ is the pairing
between the class of $E$ in $\Ext(\cO,\Lambda^{2}U_{n+1})$ and the morphism $\Lambda^{2}U_{n+1}\to \Lambda^{2}(U_{n+1}/U_{n-1})=\Omega$ and 
$\s{E,\Lambda^{2}s}$ is the pairing between the class of $E$ in $\Ext(\cO,\Lambda^{2}U_{n+1})$ and $\Lambda^{2}s:\Lambda^{2}U_{n+1}\to \Lambda^{2}M$ followed by $\Lambda^{2}M\to \Omega$. 

Let $a_n=n(n+1)(1-g)(n-\frac{1}{2})$, this is the dimension of the stack classifying extension $0\to \wedge^2 U_{n+1}\to ?\to\cO\to 0$ of $\cO_X$-modules for any fixed $(U_1,\ldots, U_{n+1})\in \cN_{0,n}$.

Let $d_{\alpha_{\cX}}$ denote the "corrected" relative dimension of $\alpha_{\cX}$, that is, $d_{\alpha_{\cX}}=a_n+\dim\cN_{0,n}+\chi(U_{n+1}^*\otimes M)$ for any $k$-points $M\in \Bun_{G_n}$ and $(U_1,\ldots, U_{n+1})\in\cN_{0,n}$.
One checks that (\ref{iso_S_Ppsi}) yields for $\cF\in \D^-(\Bun_{G_n})_!$ an isomorphism in $\D^-(\Spec k)$
$$
\Whit_{H_{n+1}}\circ F(\cF)\,\iso\,\R\Gamma_{c}(\cX,\alpha_{\cX}^{*}(\cF)\otimes \beta_{\cX}^{*}(\cL_{\psi})[d_{\alpha_{\cX}}])
$$
We will show later that Theorem~\ref{Th_1} is reduced to the following proposition. 

\begin{Pp}
\label{Pp_one}
There is a isomorphism 
$\alpha_{\cX!}(\beta_{\cX}^{*}(\cL_{\psi})[2a_n])\,\iso\,\pi_{G_n !}\epsilon _{G_n}^{*}(\cL_{\psi})$ in 
$\D^-(\Bun_{G_n})_!$.
\end{Pp}

The proposition is a consequence of the following lemmas. 
Let $\cY$ be the stack classifying  $(M,(U_{1},...,U_{n+1}),s)$ with $M\in \Bun_{G_n}$, $(U_{1},...,U_{n+1})\in \cN_{0,n}$ (i.e. $U_{k+1}/U_{k}=\Omega^{n-k}$ for $k=0,...,n$),  and $s:U_{n+1}\to M$ a morphism such that the composition 
$\Lambda^{2}U_{n+1}\vad{\Lambda^{2}s} \Lambda^{2}M \to \Omega$ coincides with $\Lambda^{2}U_{n+1}\to \Lambda^{2}(U_{n+1}/U_{n-1})=\Omega$. 

Let $\alpha_{\cY} : \cY\to \Bun_{G_n}$ be the morphism 
$(M,(U_{1},...,U_{n+1}),s)\mapsto M$. Let $\beta_{\cY} : \cY\to \AA^{1}$ be the map sending $(M,(U_{1},...,U_{n+1}),s)\in \cY$ to $\epsilon_{0,n}(U_{1},...,U_{n+1})$. 

\begin{Lm}
There is an isomorphism $\alpha_{\cX,!} \beta_{\cX}^{*}(\cL_{\psi})=\alpha_{\cY,!} \beta_{\cY}^{*}(\cL_{\psi})[-2a_n]$
 in $\D^-(\Bun_{G_n})_!$. \QED
\end{Lm}
 
 For $i\in \{1,...,n+1\}$ let $\cY_{i}$ denote the open subset of $\cY$ given by the condition that the image of $U_{i}$ by $s$ is a subbundle of $M$. One has open immersions $\cY_{n+1}\subset \cY_{n}\subset ...\subset \cY_{1}\subset \cY$. Denote by $\alpha_{\cY_{i}}:\cY_{i}\to \Bun_{G_n}$ and $\beta_{\cY_{i}}:\cY_{i}\to\AA^{1}$ the restrictions of $\alpha_{\cY}$ and $\beta_{\cY}$ to $\cY_{i}$. 
 
\begin{Lm}\label{ouverts-1}
The natural maps 
$
\alpha_{\cY_{n+1},!} \beta_{\cY_{n+1}}^{*}(\cL_{\psi}) \to 
 \alpha_{\cY_{n},!} \beta_{\cY_{n}}^{*}(\cL_{\psi}) \to ...\to 
\alpha_{\cY_{1},!} \beta_{\cY_{1}}^{*}(\cL_{\psi}) \to 
\alpha_{\cY,!} \beta_{\cY}^{*}(\cL_{\psi})
$
are isomorphisms in $\D^-(\Bun_{G_n})_!$. 
\end{Lm}
\begin{Prf}
First, one has $\cY_{n+1}=\cY_{n-1}$ thanks to the condition that the composition 
$\Lambda^{2}U_{n+1}\vad{\Lambda^{2}s} \Lambda^{2}M \to \Omega$ coincides with $\Lambda^{2}U_{n+1}\to \Lambda^{2}(U_{n+1}/U_{n-1})=\Omega$. 

Write $\cY_{0}=\cY$. Let $i\in \{1,...,n-1\}$. We are going to prove that the natural map
$$
\alpha_{\cY_{i},!} \beta_{\cY_{i}}^{*}(\cL_{\psi}) \to 
 \alpha_{\cY_{i-1},!} \beta_{\cY_{i-1}}^{*}(\cL_{\psi})
$$ 
is an isomorphism. Set $\cZ_{i}=\cY_{i-1}\backslash \cY_{i}$, let $\alpha_{\cZ_{i}} $ and $\beta_{\cZ_{i}} $ be the restrictions of $\alpha_{\cY_{i-1}} $ and $\beta_{\cY_{i-1}} $ to $\cZ_{i}$. We must prove that $\alpha_{\cZ_{i},!} \beta_{\cZ_{i}}^{*}(\cL_{\psi})=0$. 
 
  Let $\cT_{i}$ be stack classifying  
 $(M,(U_{1},U_{2},...,U_{i}),s_{i})$ with $M\in \Bun_{G_n}$, $(U_{1},U_{2},...,U_{i})\in \cN_{n-i+1,n}$, $s_{i}:U_{i}\to M$ such that the restriction of $s_{i}$ to $U_{i-1}$ is injective and its image is a subbundle of $M$, but the image of $s_{i}$ is not a subbundle of $M$ of the same rank as $U_{i}$. The map $\alpha_{\cZ_{i}}$ decomposes naturally as $\cZ_{i}\vad{\gamma_{\cZ_{i}}} \cT_{i}\vad{\alpha_{\cT_{i}}} \Bun_{G_n}$. It suffices to show that  the $*$-fibre of $\gamma_{\cZ_{i},!}\beta_{\cZ_{i}}^{*}(\cL_{\psi})$ at any closed point $(M,(U_{1},U_{2},...,U_{i}),s_{i})\in \cT_{i}$ vanishes.
 
  The fiber $\cQ$ of $\gamma_{\cZ_{i}}$ over this point is the stack classifying $((U_{1},...,U_{n+1}),s)$, where $(U_{1},...,U_{n+1})\in \cN_{0,n}$ extends $(U_{1},U_{2},...,U_{i})$, $s:U_{n+1}\to M$ extends $s_{i}$, and the composition 
$\Lambda^{2}U_{n+1}\vad{\Lambda^{2}s} \Lambda^{2}M \to \Omega$ coincides with 
$
\Lambda^{2}U_{n+1}\to \Lambda^{2}(U_{n+1}/U_{n-1})=\Omega
$.  

 Let $F$ denote the smallest subbundle of $M$ containing $s(U_{i})$, its rank is $i$ or $i-1$. Let $\cR$ be stack classifying $((W_{1},...,W_{n+1-i}),t)$ with  $(W_{1},...,W_{n+1-i})\in \cN_{0,n-i}$ and $t\in \Hom(W_{n+1-i},M/F)$. 
 There is a morphism $\rho : \cQ \to \cR$ which sends 
 $((U_{1},...,U_{n+1}),s)$ to $((U_{i+1}/U_{i},...,U_{n+1}/U_{i}),\bar s)$ where $\bar s: U_{n+1}/U_{i}\to M/F$ is the reduction of $s$. Let $\beta_{\cQ}:\cQ\to \AA^{1}$ be the restriction of $\beta_{\cZ_{i}}$ to $\cQ$. 
It suffices to show that $\rho_{!}\beta_{\cQ}^{*}(\cL_{\psi})=0$. 
 
 Pick $((W_{1},...,W_{n+1-i}),t)\in  \cR$, let $\cS$ be the fiber of $\rho$ over $((W_{1},...,W_{n+1-i}),t)$. Write $\beta_{\cS}$ for the restriction of $\beta_{\cQ}$ to $\cS$. We will show that $R\Gamma_{c}(\cS,\beta_{\cS}^{*}(\cL_{\psi}))=0$.
 
 If $F$ is of rank $i-1$ then $\cS$ identifies with the stack classifying extensions $0\to U_i/U_{i-1}\to ?\to U_{n+1}/U_i\to 0$ of $\cO_X$-modules. Since $\beta_{\cS}$ is a nontrivial character, we are done in this case.
 
 If $F$ is of rank $i$ then $\cS$ is a scheme with a free  transitive action of $\Hom(U_{n+1}/U_i, F/s(U_i))$. Under the action of $\Hom(U_{n+1}/U_i, F/s(U_i))$, $\beta_{\cS}$ changes by some character 
$$
\Hom(U_{n+1}/U_i, F/s(U_i))\to \Hom(U_{i+1}/U_i, F/s(U_i))\toup{\delta}\AA^1,
$$ 
If $D=\div(F/s(U_i))$ then $F/s(U_i)\,\iso\, \Omega^{n-i+1}(D)/\Omega^{n-i+1}$ naturally, and $\delta: \H^0(X, \Omega(D)/\Omega)\to \H^1(X,\Omega)$ is the map induced by the short exact sequence $0\to \Omega\to \Omega(D)\to \Omega(D)/\Omega\to 0$, i.e. it is the sum of the residues. Since $D>0$, $\delta$ is nontrivial, and we are done.
\end{Prf} 

\begin{Lm}
There is an isomorphism $\mu: \cY_{n+1}\to \cN_{G_n}$ such that $\pi_{G_n}\circ \mu=\alpha_{\cY_{n+1}}$ and 
$\epsilon_{G_n}\circ \mu=\beta_{\cY_{n+1}}$.\QED
\end{Lm}

 It remains to show that Proposition~\ref{Pp_one}   implies Theorem~\ref{Th_1} . 
 By the base change theorem we have 
 $$\Whit_{G_{n}}(\cF)=R\Gamma_{c}(\cN_{G_{n}},\epsilon_{G_{n}}^{*}(\cL_{\psi})\otimes \pi_{G_{n}}^{*}(\cF))[d_{N(G_{n})}-d_{G_{n}}]$$ $$=R\Gamma_{c}(\Bun_{G_{n}},\pi_{G_{n},!}\epsilon_{G_{n}}^{*}(\cL_{\psi})\otimes \cF)[d_{N(G_{n})}-d_{G_{n}}]$$ and 
 $$ R\Gamma_{c}(\cX,\alpha_{\cX}^{*}(\cF)\otimes \beta_{\cX}(\cL_{\psi})[d_{\alpha_{\cX}}])=
 R\Gamma_{c}(\Bun_{G_{n}},\alpha_{\cX,!}(\beta_{\cX}^{*}(\cL_{\psi})\otimes \cF[d_{\alpha_{\cX}}])).$$
 It remains to prove 
  $d_{\alpha_{\cX}}-2a_n=d_{N(G_n)}-d_{G_n}$. 
  This follows from 
  $d_{G_{n}}=-(1-g)n(2n+1)$, $d_{N(G_{n})}-\dim \cN_{0,n}=(1-g)(-n^{2}+n(n+1)(n-\frac{1}{2}))$, and $\chi(U_{n+1}^{*}\otimes M)=(1-g)2n^{2}(n+1)$ where $(U_{1},\dots,U_{n+1})$ and $M$ are closed points in $\cN_{0,n}$ and $\Bun_{G_{n}}$.

\subsection{From $\SO_{2n}$ to $\Sp_{2n}$}

Let $F:\D^-(\Bun_{H_n})_!\to \D^{\prec}(\Bun_{G_n})$ be the Theta functor introduced in (\cite{L1}, Definition~2). 

\begin{Th}
\label{Th_2}
The functors $\Whit_{G_n}\circ F$
and $\Whit_{H_n}$ from $\D^-(\Bun_{H_n})_!$ to $\D^-(\Spec k)$ are isomorphic. 
\end{Th}

We use the same letters as in the last paragraph (with a different meaning), as the proof is very similar.

Let $\cX$ be the stack classifying $(V,(L_{1},...,L_{n}),E,s)$ with $V\in \Bun_{H_n}$, $(L_{1},...,L_{n})\in \cN_{1,n}$ (i.e. $L_{k+1}/L_{k}=\Omega^{n-k}$ for $k=0,...,n-1$), 
an extension $0\to \Sym^{2}L_{n}\to E\to \Omega\to 0$ of $\cO_X$-modules, and a section $s:L_{n}\to V\otimes \Omega$.

Let $\alpha_{\cX} : \cX\to \Bun_{H_n}$ be the morphism 
$(V,(L_{1},...,L_{n}),E,s)\mapsto V$. 
Let $\beta_{\cX} : \cX\to \AA^{1}$ be the map sending 
$(V,(L_{1},...,L_{n}),E,s)\in \cX$ to 
$$
\epsilon_{1,n}(L_{1},...,L_{n})+\gamma(E)-\s{E,\Sym^{2}s},
$$
where $\gamma(E)$ is the pairing
between the class of $E$ in $\Ext^1(\Omega,\Sym^{2}L_{n})$ and the map $\Sym^{2}L_{n}\to \Sym^{2}(L_{n}/L_{n-1})=\Omega^{2}$; $\s{E,\Sym^{2}s}$ is the pairing between the class of $E$ in $\Ext^1(\Omega,\Sym^{2}L_{n})$ and $\Sym^{2}s:\Sym^{2}L_{n}\to \Sym^{2}V\otimes \Omega^{2}$ followed by $\Sym^{2}V\to \cO$. 

 Let $b_n=-\chi(\Omega^{-1}\otimes \Sym^{2}L_n)$ for any $k$-point $(L_1,\ldots,L_n)\in \cN_{1,n}$. Write $d_{\alpha_{\cX}}$ for the "corrected" relative dimension of $\alpha_{\cX}$, that is, 
$$
d_{\alpha_{\cX}}=\dim\cN_{1,n}+b_n+\chi(L_n^*\otimes V\otimes\Omega)
$$ 
for any $k$-points $(L_1,\ldots, L_n)\in\cN_{1,n}$ and $V\in\Bun_{H_n}$. One checks that (\ref{iso_S_Ppsi}) yields for
$\cF\in \D^-(\Bun_{H_n})_!$ an isomorphism in $\D^-(\Spec k)$
$$
\Whit_{G_n}\circ F(\cF)=\R\Gamma_{c}(\cX,\alpha_{\cX}^*(\cF)\otimes \beta_{\cX}^{*}(\cL_{\psi}))[d_{\alpha_{\cX}}] 
$$
We will derive Theorem~\ref{Th_2} from the following proposition. 

\begin{Pp}
\label{Pp_2}
There is a isomorphism 
$\alpha_{\cX,!}\beta_{\cX}^{*}(\cL_{\psi})[2b_n]\simeq \widetilde\pi_{H_n,!}\widetilde\epsilon _{H_n}^{*}(\cL_{\psi})$ in 
$\D^-(\Bun_{H_n})_!$.
\end{Pp}

Proposition~\ref{Pp_2} is reduced to the following lemmas. 
Let $\cY$ be the stack classifying $(V,(L_{1},...,L_{n}),s)$ with $V\in \Bun_{\SO_{2n}}$, $(L_{1},...,L_{n})\in \cN_{1,n}$ (i.e., $L_{k+1}/L_{k}=\Omega^{n-k}$ for $k=0,...,n-1$)
and $s:L_{n}\to V\otimes \Omega$ a morphism such that 
 the composition 
$\Sym^{2}L_{n}\vad{\Sym^{2}s} \Sym^{2}M \otimes \Omega^{2}\to \Omega^{2}$ coincides with 
$$
\Sym^{2}L_{n}
\to  \Sym^{2}(L_{n}/L_{n-1})=\Omega^{2}
$$ 
Let $\alpha_{\cY} : \cY\to \Bun_{H_n}$ be the map 
$(V,(L_{1},...,L_{n}),s)\mapsto V$. Let $\beta_{\cY} : \cY\to \AA^{1}$ be the map sending
$(V,(L_{1},...,L_{n}),s)\in \cY$ to $\epsilon_{1,n}(L_{1},...,L_{n})$.

\begin{Lm}
There is an isomorphism 
$\alpha_{\cX,!} \beta_{\cX}^{*}(\cL_{\psi})=\alpha_{\cY,!} 
\beta_{\cY}^{*}(\cL_{\psi})[-2b_n]$
 in 
$\D^-(\Bun_{H_n})_!$. \QED
\end{Lm}

 For $i\in \{1,...,n\}$ let $\cY_{i}\subset \cY$ be the open substack given by the condition that $s(L_i)\subset V\otimes \Omega$ is a subbundle of rank $i$. We have inclusions
$\cY_{n}\subset \cY_{n-1}\subset ...\subset \cY_{1}\subset \cY$. Denote by $\alpha_{\cY_{i}}:\cY_{i}\to \Bun_{H_n}$ and $\beta_{\cY_{i}}:\cY_{i}\to\AA^{1}$ the restrictions of $\alpha_{\cY}$ and $\beta_{\cY}$ to $\cY_{i}$. 

As in Lemma~\ref{ouverts-1}, one proves
\begin{Lm}
The natural maps
$\alpha_{\cY_{n},!} \beta_{\cY_{n}}^{*}(\cL_{\psi}) \to 
 \alpha_{\cY_{n-1},!} \beta_{\cY_{n-1}}^{*}(\cL_{\psi}) \to ...\to 
\alpha_{\cY_{1},!} \beta_{\cY_{1}}^{*}(\cL_{\psi}) \to 
\alpha_{\cY,!} \beta_{\cY}^{*}(\cL_{\psi}) $
are isomorphisms in $\D^-(\Bun_{H_n})_!$. \QED
\end{Lm}

\begin{Lm}
There is an isomorphism $\mu: \cY_{n}\to \widetilde \cN_{\SO_{2n}}$ such that $\widetilde\pi_{\SO_{2n}}\circ \mu=
\alpha_{\cY_{n}}$ and 
$\widetilde\epsilon_{\SO_{2n}}\circ \mu=
\beta_{\cY_{n}}$. \QED
\end{Lm}

 Theorem~\ref{Th_2} follows from Proposition~\ref{Pp_2} because  $d_{\alpha_{\cX}}-2b_n=d_{N(H_n)}-d_{H_n}$. Let us just indicate that 
 $d_{N(H_{n})}-\dim \cN_{1,n}=(1-g)n(n-1)(n-\frac{3}{2})$, 
 $\chi(L_{n}^{*}\otimes V \otimes \Omega)=(1-g)2 n^{3}$, 
 $b_{n}=(1-g)n(n+1)(n-\frac{1}{2})$ and $d_{H_{n}}=-(1-g)n(2n-1)$ where $(L_{1},\dots ,L_{n})$ and $V$ are closed points in $\cN_{1,n}$ and $\Bun_{H_{n}}$.

 \subsection{From $\GL_{n}$ to $\GL_{n+1}$}

 Let $F:\D^-(\Bun_n)_!\to \D^{\prec}(\Bun_{n+1})$ be the composition 
 of the direct  image by $\Bun_{n}\to \Bun_{n}$, $L\mapsto L^{*}$ and the theta functor $F_{n,n+1}:\D^-(\Bun_n)_!\to \D^{\prec}(\Bun_{n+1})$  introduced in (\cite{L1}, Definition~3). It is a consequence of Theorem 5 in~\cite{L1} 
 that $F$ is compatible with Hecke functors according to the morphism of dual groups $\GL_{n}\to \GL_{n+1}$, $A\mapsto \begin{pmatrix}A&0\\0&1\end{pmatrix}$.

Let us recall the definition of $F$. Denote $\cW$  be the classifying stack of $(L,U,s)$ with $L\in \Bun_{n}$, $U\in \Bun_{n+1}$ and $s:L\to U$ a morphism. We have 
$(h_{n},h_{n+1}):\cW\to \Bun_{n}\times \Bun_{n+1}, (L,U,s)\mapsto (L,U)$. Then for $\cF\in \D^-(\Bun_{n})_!$, 
$$
F(\cF)=h_{n+1,!}((h_{n}^{*}\cF)[\dim\Bun_{n+1}+\chi(L^{*}\otimes U)]), 
$$
where $\chi(L^{*}\otimes U)$ is considered as a locally constant function on $\Bun_{n}\times \Bun_{n+1}$.

\begin{Th}
\label{Th_3}
The functors $\Whit_{\GL_{n+1}}\circ F$
and $\Whit_{\GL_{n}}$ from $\D^-(\Bun_{n})_!$ to $\D-(\Spec k)$ are isomorphic. 
\end{Th}

Let $\cX$ be the stack classifying 
$L\in \Bun_{n}$, 
$(U_{1},\dots ,U_{n+1})\in \cN_{0,n}$, and $s:L\to U_{n+1}$ a morphism. We have $\alpha_{\cX}:\cX \to \Bun_{n}$ and 
$\beta_{\cX}:\cX \to \mathbb A^{1}$ which send $(L,(U_{1},\dots U_{n+1}),s)$ to $L$ and $\epsilon_{0,n}(U_{1},\dots U_{n+1})$. 

We have 
$$\Whit_{\GL_{n+1}}\circ F(\cF)=R\Gamma_{c}(\Bun_{n},\cF\otimes \alpha_{\cX,!}\beta^{*}_{\cX}(\cL_{\psi})[\dim\cN_{0,n}+\chi(L^{*}\otimes U_{n+1})])\ \ \text{and}$$
$$\Whit_{\GL_{n}}(\cF)=R\Gamma_{c}(\Bun_{n},\cF\otimes
(\pi_{0,n-1})_{!}\epsilon_{0,n-1}^{*}(\cL_{\psi})[\dim\cN_{0,n-1}-\dim\Bun_{n}]).$$

 For $i\in \{0,\dots ,n\}$ denote by $\cX_{i}$ the open substack of $\cX$ classifying $(L,(U_{1},\dots U_{n+1}),s)$ such that the composition $L\vad{s} U_{n+1}\to U_{n+1}/U_{n+1-i}$ is surjective. We have $\cX=\cX_{0}\supset \cX_{1} \supset \dots \supset \cX_{n}$ and we have an isomorphism  
 $\cN_{0,n-1}\to \cX_{n}$ which sends $(E_{1}, \dots ,E_{n})$ to 
 $(E_{n}, (\Omega^{n}, \Omega^{n}\oplus E_{1},\dots ,\Omega^{n}\oplus E_{n}),(0,\mathrm{Id}))$ with
  $(0,\mathrm{Id}):E_{n}\to \Omega^{n}\oplus E_{n}$ the obvious inclusion. 
 
 It is easy to compute that for $L=E_{n}$ with $(E_{1}, \dots ,E_{n})\in \cN_{0,n-1}$ and 
 $(U_{1},\dots ,U_{n+1})\in \cN_{0,n}$ 
  we have 
 $\dim\cN_{0,n}+\chi(L^{*}\otimes U_{n+1})
 =\dim\cN_{0,n-1}-\dim\Bun_{n}$. 
 
 Therefore we are reduced to the following lemma. 
 We denote by $\alpha_{\cX_{i}}:\cX_{i} \to \Bun_{n}$ and 
$\beta_{\cX_{i}}:\cX_{i} \to \mathbb A^{1}$ the restrictions of $\alpha_{\cX}$ and $\beta_{\cX}$ to $\cX_{i}$.

 \begin{Lm}\label{ouverts-GL}
The natural maps 
$
\alpha_{\cX_{n},!} \beta_{\cX_{n}}^{*}(\cL_{\psi}) \to 
 \alpha_{\cX_{n-1},!} \beta_{\cX_{n-1}}^{*}(\cL_{\psi}) \to ...\to 
\alpha_{\cX_{1},!} \beta_{\cX_{1}}^{*}(\cL_{\psi}) \to 
\alpha_{\cX,!} \beta_{\cX}^{*}(\cL_{\psi})
$
are isomorphisms in $\D^-(\Bun_{n})_!$. 
\end{Lm}
\begin{Prf}
 We recall that $\cX=\cX_{0}$. 
  Let $i\in \{1,...,n\}$. We are going to prove that the natural map
$$
\alpha_{\cX_{i},!} \beta_{\cX_{i}}^{*}(\cL_{\psi}) \to 
 \alpha_{\cX_{i-1},!} \beta_{\cX_{i-1}}^{*}(\cL_{\psi})
$$ 
is an isomorphism. Set $\cZ_{i}=\cX_{i-1}\backslash \cX_{i}$, let $\alpha_{\cZ_{i}} $ and $\beta_{\cZ_{i}} $ be the restrictions of $\alpha_{\cX_{i-1}} $ and $\beta_{\cX_{i-1}} $ to $\cZ_{i}$. We must prove that $\alpha_{\cZ_{i},!} \beta_{\cZ_{i}}^{*}(\cL_{\psi})=0$. 
 
  Let $\cT_{i}$ be stack classifying  
 $(L,(V_{1},V_{2},...,V_{i}),t)$ with $L\in \Bun_n$, $(V_{1},V_{2},...,V_{i})\in \cN_{0,i-1}$, $t:L\to V_{i} $ such that the
 composition $L\vad{t} V_{i}\to V_{i}/V_{1}$ is surjective  but $t$ is not surjective. The map $\alpha_{\cZ_{i}}$ decomposes naturally as $\cZ_{i}\vad{\gamma_{\cZ_{i}}} \cT_{i}\vad{\alpha_{\cT_{i}}} \Bun_{n}$
 where $\gamma_{\cZ_{i}}(L,(U_{1},\dots U_{n+1}),s)=
 (L,(U_{n+2-i}/U_{n+1-i},\dots, U_{n+1}/U_{n+1-i}),\bar s)$ and 
 $\alpha_{\cT_{i}}(L,(U_{1},\dots U_{n+1}),s)=L$. 
  It suffices to show that  the $*$-fibre of $\gamma_{\cZ_{i},!}\beta_{\cZ_{i}}^{*}(\cL_{\psi})$ at any closed point $(L,(V_{1},V_{2},...,V_{i}),t)\in \cT_{i}$ vanishes.
 
  Let us choose a closed point $(L,(V_{1},V_{2},...,V_{i}),t)\in \cT_{i}$  and define    $L''=\Ker t$  and $L'$ the kernel of the composition $L\vad{t} V_{i}\to V_{i}/V_{1}$. Then $L'$ is a subbundle of $L$ of rank $n+1-i$ and $L''$ is a subbundle of $L$ of rank $n+1-i$ or $n-i$.

  The fiber $\cQ$ of $\gamma_{\cZ_{i}}$ over this closed point is the stack classifying $((U_{1},...,U_{n+1}),s)$, with an isomorphism between $U_{n+1}/U_{n+1-i}$ and $V_{i}$ sending $U_{j+n+1-i}/U_{n+1-i}$ to $V_{j}$ for any $j\in \{0,\dots ,i\}$ and 
  $s:L\to U_{n+1}$ such that the composition 
  $L\vad{s} U_{n+1} \to U_{n+1}/U_{n+1-i}\simeq V_{i}$ is $t$. 
  Let $\cR$ be stack classifying $((U_{1},\dots,U_{n+1-i}),s_{i})$ with  $(U_{1},\dots,U_{n+1-i})\in \cN_{i,n}$ and $s_{i}\in \Hom(L'',U_{n+1-i})$. 
 There is a morphism $\rho : \cQ \to \cR$ which sends 
 $((U_{1},...,U_{n+1}),s)$ to $((U_{1},\dots,U_{n+1-i}),s_{i})$ where $s_{i}$ is the restriction of $s$ to $L''$. 
 Let $\beta_{\cQ}:\cQ\to \AA^{1}$ be the restriction of $\beta_{\cZ_{i}}$ to $\cQ$. 
It suffices to show that $\rho_{!}\beta_{\cQ}^{*}(\cL_{\psi})=0$. 
 
 Pick $((U_{1},\dots,U_{n+1-i}),s_{i})\in  \cR$, let $\cS$ be the fiber of $\rho$ over $((U_{1},\dots,U_{n+1-i}),s_{i})$. Write $\beta_{\cS}$ for the restriction of $\beta_{\cQ}$ to $\cS$. We will show that $R\Gamma_{c}(\cS,\beta_{\cS}^{*}(\cL_{\psi}))=0$.
 
 If $L'=L''$ we have an exact sequence 
 $0\to L/L''\to U_{n+1}/U_{n+1-i}\to U_{n+2-i}/U_{n+1-i}\to 0$ and  
  $\cS$ identifies with the stack classifying extensions $0\to U_{n+1-i}\to ?\to U_{n+2-i}/U_{n+1-i}\to 0$ of $\cO_X$-modules. Since $\beta_{\cS}$ is a nontrivial character, we are done in this case.
 
 If $L'/L''$ is a line bundle then $\cS$ is a scheme with a  free transitive action of 
 the $H^{0}$ of the cone of the morphism of complexes of $k$-vector spaces 
 $$\mathrm{RHom}(U_{n+1}/U_{n+1-i},U_{n+1-i})\to 
 \mathrm{RHom}(L/L'',U_{n+1-i})$$ 
 which is also the cone of the morphism of complexes 
 $$\mathrm{RHom}(U_{n+2-i}/U_{n+1-i},U_{n+1-i})\to 
 \mathrm{RHom}(L'/L'',U_{n+1-i}) $$
 and whose cohomology is concentrated in degree $0$.  
The last morphism of complexes comes from the  non zero morphism 
$L'/L''\to U_{n+2-i}/U_{n+1-i}=\Omega^{i-1}$ which identifies 
$L'/L''$ to $\Omega^{i-1}(-D)$ for some effective non zero divisor $D$. Therefore the $H^{0}$ of this cone is equal to 
$$H^{0}(U_{n+1-i}\otimes \Omega^{1-i}(D)/U_{n+1-i}\otimes \Omega^{1-i})$$
and $\beta_{\cS}^{*}(\cL_{\psi})$ transforms under this action through the character 
$$H^{0}(U_{n+1-i}\otimes \Omega^{1-i}(D)/U_{n+1-i}\otimes \Omega^{1-i})\to $$ $$
H^{0}((U_{n+1-i}/U_{n-i})\otimes \Omega^{1-i}(D)/(U_{n+1-i}/U_{n-i})\otimes \Omega^{1-i})
=H^{0}(\Omega(D)/\Omega)\vad{\sigma} \AA^{1}$$
where $\sigma$ is the sum of the residues. 
Since $D$ is non zero, $\sigma$ is a non zero character and 
we are done.  
 \end{Prf}

\end{document}